\documentclass[12pt,a4paper,twoside]{article}

\newcommand{\pageformat}[6]{\setlength{\hoffset}{-1in}
                  \setlength{\voffset}{-1in}
                  \addtolength{\hoffset}{#5}
                            \addtolength{\voffset}{#6}
                            \setlength{\oddsidemargin}{#1}
                            \setlength{\evensidemargin}{#2}
                            \setlength{\textwidth}{\paperwidth}
                  \addtolength{\textwidth}{-\oddsidemargin}
                  \addtolength{\textwidth}{-\evensidemargin}
                  \addtolength{\textwidth}{-\marginparsep}
                  \addtolength{\textwidth}{-\marginparwidth}
                            \setlength{\topmargin}{#3}
                            \setlength{\textheight}{\paperheight}
                  \addtolength{\textheight}{-\topmargin}
                  \addtolength{\textheight}{-\headheight}
                  \addtolength{\textheight}{-\headsep}
                  \addtolength{\textheight}{-\footskip}
                  \addtolength{\textheight}{-#4}}
\pageformat{2cm}{3cm}{25mm}{25mm}{1pt}{0pt}

\usepackage{ifthen}
\newboolean{article}
    \setboolean{article}{true}
\newboolean{report}
\newboolean{book}
\newboolean{letter}
\newboolean{german}
\newboolean{italian}
\newboolean{nobaselinestretch}
\newboolean{nosectionappendix}
\newboolean{oldtoc}
\newboolean{nosectionequn}
\newboolean{notheorem}

\ifthenelse{\boolean{german}}{
    \usepackage{german}}{}

\usepackage[latin1]{inputenc}

\usepackage{amsmath}
\usepackage{amssymb}
\usepackage[mathscr]{eucal}

\ifthenelse{\boolean{notheorem}}{}{
    \usepackage{theorem}}

\ifthenelse{\boolean{nobaselinestretch}}{}{
    \renewcommand{\baselinestretch}{1.25}}

\newenvironment{env}[2]{\begin{#1}#2\end{#1}}{}
    \newcommand{\beq}[1]{\begin{env}{equation}{#1}}
    \newcommand{\beqn}[1]{\begin{env}{equation*}{#1}}
    \newcommand{\bal}[1]{\begin{env}{align}{#1}}
    \newcommand{\baln}[1]{\begin{env}{align*}{#1}}
    \newcommand{\bga}[1]{\begin{env}{gather}{#1}}
    \newcommand{\bgan}[1]{\begin{env}{gather*}{#1}}
    \newcommand{\bflal}[1]{\begin{env}{flalign}{#1}}
    \newcommand{\bflaln}[1]{\begin{env}{flalign*}{#1}}
    \newcommand{\bmu}[1]{\begin{env}{multline}{#1}}
    \newcommand{\bmun}[1]{\begin{env}{multline*}{#1}}
    \newcommand{\bsp}[1]{\begin{env}{split}{#1}}

    \newcommand{\eeq}{\end{env}}
    \newcommand{\eeqn}{\end{env}}
    \newcommand{\eal}{\end{env}}
    \newcommand{\ealn}{\end{env}}
    \newcommand{\ega}{\end{env}}
    \newcommand{\egan}{\end{env}}
    \newcommand{\eflal}{\end{env}}
    \newcommand{\eflaln}{\end{env}}
    \newcommand{\emu}{\end{env}}
    \newcommand{\emun}{\end{env}}
    \newcommand{\esp}{\end{env}}

\newcommand{\lf}{\vspace{2ex}}

\renewcommand{\bf}[1]{\textbf{#1}}
\renewcommand{\it}[1]{\textit{#1}}

\renewcommand{\sf}[1]{\textsf{#1}}

\renewcommand{\tt}[1]{\texttt{#1}}
\newcommand{\hl}[1]{\bf{\it{#1}}}
\newcommand{\mrm}[1]{\mathrm{#1}}
\newcommand{\mbf}[1]{\mathbf{#1}}
\newcommand{\msf}[1]{\text{\small$\sf{#1}$}}

\newcommand{\cmc}[1]{\mathcal{#1}}
\newcommand{\eus}[1]{\mathscr{#1}}
\newcommand{\euf}[1]{\mathfrak{#1}}
\newcommand{\bb}[1]{\mathbb{#1}}

\newcommand{\nbd}[1]{$#1$\nobreakdash--}
\newcommand{\ol}[1]{\overline{#1}}

\newcommand{\ve}{\varepsilon}
\newcommand{\vt}{\vartheta}

\newcommand{\vp}{\varphi}
\newcommand{\om}{\omega}

\newcommand{\norm}[1]{\left\lVert#1\right\rVert}

\newcommand{\Babs}[1]{\Bigl\lvert#1\Bigr\rvert}
\newcommand{\Bnorm}[1]{\Bigl\lVert#1\Bigr\rVert}
\newcommand{\snorm}[1]{\norm{\smash{#1}}}

\newcommand{\bfam}[1]{\bigl(#1\bigr)}

\newcommand{\AB}[1]{\langle#1\rangle}

\newcommand{\CB}[1]{\{#1\}}
\newcommand{\bCB}[1]{\bigl\{#1\bigr\}}
\newcommand{\BCB}[1]{\Bigl\{#1\Bigr\}}
\newcommand{\SB}[1]{[#1]}

\newcommand{\RO}[1]{[#1)}

\newcommand{\Matrix}[1]{\begin{pmatrix}#1\end{pmatrix}}

\newcommand{\set}[2][]{
    \ifthenelse{\equal{#1}{}}{
        \CB{#2}}{
        \CB{#1~|~#2}}}
\newcommand{\bset}[2][]{
    \ifthenelse{\equal{#1}{}}{
        \bCB{#2}}{
        \bCB{#1~|~#2}}}
\newcommand{\Bset}[2][]{
    \ifthenelse{\equal{#1}{}}{
        \BCB{#2}}{
        \BCB{#1~\big|~#2}}}
\newcommand{\zero}{\CB{0}}

\DeclareMathOperator{\ls}{\normalfont\msf{span}}
\DeclareMathOperator{\cls}{\ol{\ls}}

\DeclareMathOperator{\id}{\normalfont\msf{id}}

\renewcommand{\Re}{\operatorname{\msf{Re}}}
\renewcommand{\Im}{\operatorname{\msf{Im}}}

\newcommand{\C}{\bb{C}}

\newcommand{\N}{\bb{N}}

\newcommand{\R}{\bb{R}}

\newcommand{\cA}{\cmc{A}}
\newcommand{\cB}{\cmc{B}}

\newcommand{\cF}{\cmc{F}}

\newcommand{\sB}{\eus{B}}

\newcommand{\sK}{\eus{K}}

\newcommand{\sS}{\eus{S}}

\newcommand{\eH}{\euf{H}}

\newcommand{\eK}{\euf{K}}

\newcommand{\eT}{\euf{T}}

\newcommand{\U}{\mbf{1}}

\newcommand{\G}{\Gamma}

\newcommand{\DG}{{\mrm{I}\hspace{-0.3ex}\G}}

\newcommand{\I}{{I\!\!\!\;I}}

\newcommand{\s}{\text{\scriptsize$\sS$}}

\newcommand{\0}{\mbf{0}}

\ifthenelse{\boolean{nosectionequn}}{}{
    \numberwithin{equation}{section}
    }

\ifthenelse{\boolean{article}\or\boolean{letter}\or\boolean{nosectionequn}}{
    \setboolean{nosectionappendix}{true}}{}
\ifthenelse{\boolean{nosectionappendix}}{}{
    \renewcommand{\appendix}{
        \chapter*{\appendixname}
        \addcontentsline{toc}{chapter}{\appendixname}
        \renewcommand{\thesection}{\Alph{section}}
        \setcounter{section}{0}}}
   
\ifthenelse{\boolean{report}\or\boolean{book}}{
    }{}

\ifthenelse{\boolean{notheorem}}{}{
        \newcommand{\mnname}{Mathematical note.}
        \newcommand{\enname}{End of the note.}
        \newcommand{\definame}{Definition.}
        \newcommand{\propname}{Proposition.}
        \newcommand{\lemname}{Lemma.}
        \newcommand{\exname}{Example.}
        \newcommand{\exername}{Exercise.}
        \newcommand{\remname}{Remark.}
        \newcommand{\obname}{Observation.}
        \newcommand{\thmname}{Theorem.}
        \newcommand{\corname}{Corollary.}
        \newcommand{\proofname}{Proof.}
    \ifthenelse{\boolean{german}}{
        \renewcommand{\mnname}{Mathematische Notiz.}
        \renewcommand{\enname}{Ende der Notiz.}
        \renewcommand{\exname}{Beispiel.}
        \renewcommand{\exername}{Übung.}
        \renewcommand{\remname}{Bemerkung.}
        \renewcommand{\obname}{Beobachtung.}
        \renewcommand{\thmname}{Satz.}
        \renewcommand{\corname}{Korollar.}
        \renewcommand{\proofname}{Beweis.}}{}
    \ifthenelse{\boolean{italian}}{
        \renewcommand{\mnname}{Nota matematica.}
        \renewcommand{\enname}{Fina della nota.}
        \renewcommand{\definame}{Definizione.}
        \renewcommand{\propname}{Proposizione.}
        \renewcommand{\exname}{Esempio.}
        \renewcommand{\exername}{Esercizio.}
        \renewcommand{\remname}{Nota.}
        \renewcommand{\obname}{Osservazione.}
        \renewcommand{\thmname}{Teorema.}
        \renewcommand{\corname}{Corollario.}
        \renewcommand{\proofname}{Dimostrazione.}

       \renewcommand{\appendixname}{Appendice}

       }{}
    \theoremheaderfont{\normalfont\bfseries}
    \theoremstyle{change}
        \theorembodyfont{\rmfamily}
            \newtheorem{emp}{}[section]
                \newcommand{\bemp}[1][]{
                    \begin{emp}\hskip-\labelsep\bf{#1}\hskip\labelsep}
                \newcommand{\eemp}{\end{emp}}
\newtheorem{itemp}[emp]{}
                \newcommand{\bitemp}[1][]{
                    \begin{itemp}\hskip-\labelsep\bf{#1}\hskip\labelsep\normalfont\itshape}
                \newcommand{\eitemp}{\end{itemp}}
            \newtheorem{mn}[emp]{\mnname}
                \newcommand{\bnm}{\begin{mn}~\begin{quotation}\renewcommand{\baselinestretch}{1}\small\noindent\ignorespaces}
                \newcommand{\enm}{\end{quotation}\hfill\bf{\enname}\end{mn}}
            \newtheorem{ex}[emp]{\exname}
                \newcommand{\bex}{\begin{ex}}
                \newcommand{\eex}{\end{ex}}
            \newtheorem{exer}[emp]{\exername}
                \newcommand{\bexer}{\begin{exer}}
                \newcommand{\eexer}{\end{exer}}
            \newtheorem{defi}[emp]{\definame}
                \newcommand{\bdefi}{\begin{defi}}
                \newcommand{\edefi}{\end{defi}}
            \newtheorem{rem}[emp]{\remname}
                \newcommand{\brem}{\begin{rem}}
                \newcommand{\erem}{\end{rem}}
            \newtheorem{ob}[emp]{\obname}
                \newcommand{\bob}{\begin{ob}}
                \newcommand{\eob}{\end{ob}}
        \theorembodyfont{\normalfont\itshape}
            \newtheorem{thm}[emp]{\thmname}
                \newcommand{\bthm}{\begin{thm}}
                \newcommand{\ethm}{\end{thm}}
            \newtheorem{prop}[emp]{\propname}
                \newcommand{\bprop}{\begin{prop}}
                \newcommand{\eprop}{\end{prop}}
            \newtheorem{cor}[emp]{\corname}
                \newcommand{\bcor}{\begin{cor}}
                \newcommand{\ecor}{\end{cor}}
            \newtheorem{lem}[emp]{\lemname}
                \newcommand{\blem}{\begin{lem}}
                \newcommand{\elem}{\end{lem}}
\newenvironment{empn}[1]{\lf\noindent\bf{#1}\ignorespaces\hskip\labelsep}{\lf}
		\newcommand{\bempn}[1]{\begin{empn}{#1}}
		\newcommand{\eempn}{\end{empn}}
		\newcommand{\bitempn}[1]{\begin{empn}{#1}\normalfont\itshape}
		\newcommand{\eitempn}{\end{empn}}
                \newcommand{\bnmn}{\begin{empn}{\mnname}~\begin{quotation}\renewcommand{\baselinestretch}{1}\small\noindent\ignorespaces}
                \newcommand{\enmn}{\end{quotation}\hfill\bf{\enname}\end{empn}}
		\newcommand{\bexn}{\begin{empn}{\exname}}
		\newcommand{\eexn}{\end{empn}}
		\newcommand{\bexern}{\begin{empn}{\exername}}
		\newcommand{\eexern}{\end{empn}}
		\newcommand{\bdefin}{\begin{empn}{\definame}}
		\newcommand{\edefin}{\end{empn}}
		\newcommand{\bremn}{\begin{empn}{\remname}}
		\newcommand{\eremn}{\end{empn}}
		\newcommand{\bobn}{\begin{empn}{\obname}}
		\newcommand{\eobn}{\end{empn}}

\newcommand{\qedsymbol}{~\rule[-0.35mm]{2mm}{2mm}}
    \newcounter{proof}[emp]
    \newenvironment{Proof}[1]{
        \vspace{1ex}
        \renewcommand{\item}[1][\stepcounter{proof}(\roman{proof})]%
            {##1\hskip\labelsep}
        \noindent\textsc{#1\hskip\labelsep}}{
        \nolinebreak\qedsymbol}
    \newcommand{\proof}[1][\proofname]{
        \begin{Proof}{#1}\ignorespaces}
    \newcommand{\qed}{\end{Proof}}
    \newcommand{\noqed}{
        \renewcommand{\qedsymbol}{}
        \end{Proof}}}
    \ifthenelse{\boolean{italian}}{
        \renewcommand{\proofname}{Dimostrazione.}}{}

\usepackage[varg]{txfonts}

\usepackage[hypertex]{hyperref}

\begin{document}

\bibliographystyle{amsalpha}

\title{Subsystems of Fock Need Not Be Fock:\\Spatial CP-Semigroups\thanks{This work has been supported by an RiP-Program at Mathematisches Forschungsinstitut Oberwolfach. BVRB is supported by the Department of Science and Technology (India) under the Swarnajayanthi Fellowship Project. MS is supported by research funds of the Dipartimento S.E.G.e S.\ of University of Molise and of the Italian MUR (PRIN 2005 and 2007).}}

\author{B.V.Rajarama Bhat, Volkmar Liebscher, and Michael Skeide}

\date{April 2008; revised July 2009}

\maketitle

\begin{abstract}
\noindent
We show that a product subsystem of a time ordered system (that is, a product system of time ordered Fock modules), though type I, need not be isomorphic to a time ordered product system. In that way, we answer an open problem in the classification of CP-semigroups by product systems. We define spatial strongly continuous CP-semigroups on a unital \nbd{C^*}algebra and characterize them as those that have a Christensen-Evans generator.
\end{abstract}

\section{Introduction}

Quantum dynamics deals with \it{\nbd{E_0}semigroups} (that is, semigroups of unital endomorphisms of a $C^*$ or von Neumann algebra) or, more generally, \it{CP-semigroups} (that is, semigroups of completely positive maps), and tries to classify them in terms of product systems.

Let $H$ denote a Hilbert space. Product systems of Hilbert spaces (\hl{Arveson systems} for short), \nbd{E_0}semigroups on $\sB(H)$, and CP-semigroups on $\sB(H)$ are intimately related among each other. For all three there is a notion of \it{spatiality}; spatial Arveson systems (Arveson \cite{Arv89}), spatial \nbd{E_0}semigroups (Powers \cite{Pow87}), and spatial CP-semigroups (Arveson \cite{Arv97a}). It is well known that these properties are equivalent, when suitably formulated.

There is a similar relation between product systems of correspondences (that is, Hilbert bimodules) over a \nbd{C^*}algebra $\cB$, \nbd{E_0}semigroups acting on the algebra $\sB^a(E)$ of adjointable operators on a Hilbert \nbd{\cB}module $E$, and CP-semigroups on $\cB$. Generalizing Powers definition to $\sB^a(E)$, spatial \nbd{E_0}semigroups correspond to \it{spatial} product systems (Skeide \cite{Ske06d}). A similar analysis for \it{spatial} CP-semigroups and their product systems is missing.

In these notes we define spatial CP-semigroups on a unital \nbd{C^*}algebra $\cB$ in analogy to Arveson's definition: A CP-semigroup $T=\bfam{T_t}_{t\in\R_+}$ on $\cB$ is \hl{spatial}, if there exists a continuous semigroup $c=\bfam{c_t}_{t\in\R_+}$ in $\cB$ such that for all $t\in\R_+$ the map $T_t-c_t^*\bullet c_t$ on $\cB$ is completely positive; see Section \ref{spatialCP}. It turns out that CP-semigroups are spatial if and only if their associated product system, the \it{GNS-system}, embeds into a spatial one (Theorem \ref{CPspatPSthm}). By a counter example in the central Section \ref{counter} we show that this result cannot be improved. Another important result is that a CP-semigroup is spatial if and only if its generator has \it{Christensen-Evans form} (Corollary \ref{spatCEcor}). Actually, our counter example in Section \ref{counter} shows more: Even if a subsystem of a product system of \it{time ordered Fock modules} is generated by a single continuous \it{unit}, then it need not be isomorphic to a product system of time ordered Fock modules. This is in sharp contrast with the Hilbert space case. Indeed, well-known results assert that a subsystem of a \it{type I} product system of Hilbert spaces is type I, too, and that a type I system is isomorphic to a product systems of symmetric Fock spaces. Even worse, in Theorem \ref{tInspatthm} we show that a type I product system need not be contained in a product system of time ordered Fock modules. These are important results in the classification of product systems.

\lf
We would like to thank the referee for having encouraged this revision, making these notes more self-contained. This also gave us the occasion to illustrate in some places how, meanwhile, the notions and results of these notes have been applied and generalized into several directions. The most important of these results is that spatial CP-semigroups on a von Neumann algebra (where $c$ is only required to be strongly continuous) are precisely those that have spatial product systems of \it{von Neumann correspondences}; see Skeide \cite{Ske08p1}. This means that most peculiarities disappear when passing to von Neumann algebras; see Section \ref{rem}.

In Section \ref{prel} we start right away with the definitions and known facts that regard our notes, that is, that regard the \nbd{C^*}case. Apart from the discussion of the von Neumann case, more remarks regarding both the general context and newer results are postponed to Section \ref{rem}, too.

\section{Preliminaries}\label{prel}

\bdefi
Let $\cB$ denote a \nbd{C^*}algebra. A \hl{product system} is a family $E^\odot=\bfam{E_t}_{t\in\R_+}$ of correspondences $E_t$ over $\cB$ together with a family of bilinear unitaries
\beqn{
u_{s,t}
\colon
E_s\odot E_t\rightarrow E_{s+t}
}\eeqn
such that the \hl{product} $x_sy_t:=u_{s,t}(x_s\odot y_t)$ is associative. (By $\odot$ we denote the \it{internal tensor product} over $\cB$.) We always require that $E_0=\cB$ and that $u_{0,t}$ and $u_{s,0}$ are left and right multiplication, respectively, with elements in $\cB=E_0$.
\edefi

Apart from the marginal conditions for $u_{0,t}$ and $u_{s,0}$, this is Bhat and Skeide \cite[Definition 4.7]{BhSk00}. Note that this definition ignores possible technical conditions like continuity, or measurability with respect to time $t\in\R_+$, or (quasi) triviality of the bundle structure. If we wish to emphasize that there are no technical conditions, we sometimes refer to $E^\odot$ as \hl{algebraic} product system. We do \bf{not} follow some authors who would say \it{discrete} product system. (In our thinking \it{discrete} refers to product systems indexed by $\N_0$.) Recall that \hl{Arveson system} refers to the case $\cB=\C$, where we also use the usual $\otimes$ sign.

Product systems are classified, in a first step, by their supply of units.

\bemp[{\protect\cite[Definition 4.7]{BhSk00}}.~]
A \hl{unit} for a product system $E^\odot$ is a family $\xi^\odot=\bfam{\xi_t}_{t\in\R_+}$ of elements $\xi_t\in E_t$ that is multiplicative in the sense that $\xi_s\xi_t=\xi_{s+t}$, and where $\xi_0=\U\in\cB$.

A unit $\xi^\odot$ is \hl{unital}, if it consists of \hl{unit vectors}, that is, if $\AB{\xi_t,\xi_t}=\U$ for all $t\in\R_+$.
\eemp

In order to find satisfactory classification results, the units have to satisfy technical conditions. The best correspondence to the case of Arveson systems occurs, if we require the units to be \it{continuous}. Directly from the definition of the tensor product, it follows that for every pair of units $\xi^\odot,\xi'^\odot$ the bounded maps $\eT^{\xi,\xi'}_t=\AB{\xi_t,\bullet\xi'_t}$ form a semigroup $\eT^{\xi,\xi'}=\bfam{\eT^{\xi,\xi'}_t}_{t\in\R_+}$ on $\cB$. The following definition is from Skeide \cite{Ske03b}.

\bdefi
As set $S$ of units for a product system $E^\odot$ is \hl{continuous}, if $\eT^{\xi,\xi'}$ is uniformly continuous for all $\xi^\odot,\xi'^\odot\in S$. A single unit $\xi^\odot$ is \hl{continuous}, if $\bCB{\xi^\odot}$ is continuous.

A product system is \hl{type I}, if it is generated by a continuous set of units $S$. It is \hl{type II} if it has a continuous unit, but is not type I. It is \hl{type III}, if it has no continuous units at all.
\edefi

The property of a set $S$ of units to be continuous is intrinsic to the family $\eT^{\xi,\xi'}$ of semigroups. It is not related to a whatsoever continuous structure $E^\odot$ might have. However, if $E^\odot$ has a continuous structure, then one would insist that $S$ consists of sections that are continuous also with respect to that continuous structure. 

\it{Continuous product systems} have been defined in Skeide \cite[Definition 7.1]{Ske03b}. We do not reproduce that definition here, since we are only interested in the case with units. We just mention that a part of the definition is a distinguished set of \hl{continuous sections} $CS(E^\odot)$ that turns $E^\odot$ into a \it{continuous field} of Banach spaces; see Dixmier \cite{Dix77}. And as soon as there is a continuous unital unit among the continuous sections, the continuous structure of $E^\odot$ may be recovered from that unit. The following facts are a slight reformulation of \cite[Theorem 7.5]{Ske03b}, taking into account also (the proof of) \cite[Theorem 10.2]{BhSk00} and \cite[Lemma 4.4.11]{BBLS04}.

\bthm\label{unitCPSthm}
\begin{enumerate}
\item
If $E^\odot$ is a product system and $\xi^\odot$ a continuous unit, then there is at most one continuous structure $CS(E^\odot)$ on $E^\odot$ such that $\xi^\odot\in CS(E^\odot)$.

\item
If $E^\odot$ is generated by a continuous set of units $S$, then it possesses a unique continuous structure $CS(E^\odot)\supset S$.
\end{enumerate}
\ethm

\noindent
The quoted results are under the hypothesis that at least one of the involved units is unital. But, continuous units may be normalized suitably; see the second construction in Liebscher and Skeide \cite[Example 4.2]{LiSk08}. More precisely, whatever the continuous set of units $S$ is, we may assume that it contains at least one unital unit without changing the product subsystem generated by $S$. On the other hand, a look at the involved proofs shows that Theorem \ref{unitCPSthm} remains true, if we speak about \it{strongly continuous} sets of units $S$ in the sense that all semigroups $\eT^{\xi,\xi'}$ are strongly continuous, provided there is at least one unital unit $\xi^\odot\in S$. Strongly continuous units may not necessarily be normalized.

We see that for non-type III product systems it is enough to fix a single continuous unit, in order to fix, if it exists, the whole continuous structure. However, it is not known if there are nonisomorphic continuous structures on the same algebraic product system.\footnote{For Arveson systems, Liebscher \cite[Corollary 7.16]{Lie09} states that the algebraic structure determines the measurable structure up to isomorphism. We do not know (not even for Arveson systems!), if this remains true for the continuous structure of an Arveson system.} In the sequel, we always think of continuous product systems and require $S\subset CS(E^\odot)$.

If $E^\odot$ is a continuous product system, then the map $t\mapsto\AB{x_t,by_t}$ is continuous for all $b\in\cB$ and $x,y\in CS(E^\odot)$. The following theorem illustrates how strong the assumption of existence of a continuous unit among the continuous sections really is.

\bitemp[{\protect\cite[Theorem 7.7]{Ske03b}}.~]\label{ucCPSthm}
If $E^\odot$ is a continuous product system with a continuous unit $\xi^\odot\in CS(E^\odot)$, then $E^\odot$ is \hl{uniformly continuous} in the sense that the map $t\longmapsto\AB{x_t,\bullet y_t}$ is uniformly continuous for all $x,y\in CS(E^\odot)$.
\eitemp

We now come to \it{spatial} product systems.

\bemp[Definition {\protect\cite[Section 2]{Ske06d}}.~]
A unit $\om^\odot$ is \hl{central}, if $b\om_t=\om_tb$ for all $b\in\cB,t\in\R_+$. A product system $E^\odot$ is \hl{spatial}, if it admits a central unital unit.\footnote{Actually, like in \cite{Ske06d} one should speak of pairs $(E^\odot,\om^\odot)$, because the \it{spatial structure} may depend on the choice of $\om^\odot$; see Tsirelson \cite{Tsi08}.} If $E^\odot$ is continuous, then we require that $\om^\odot\in CS(E^\odot)$.
\eemp

\bob\label{spatUCPSob}
If $E^\odot$ is spatial, then the reference unit $\om^\odot$ is, clearly, continuous. ($\eT^{\om,\om}$ is the trivial semigroup.) By Theorem \ref{ucCPSthm}, a spatial continuous product system $E^\odot$ is uniformly continuous.
\eob

We now explain the differences in the classification of Arveson systems and general (continuous) product systems. Most statements are based on the two fundamental Examples \ref{CPGNSex} and \ref{E0BPSex}, and on time ordered product systems as discussed in Section \ref{counter}.

\bex\label{CPGNSex}
A \hl{CP-semigroup} is a semigroup $T=\bfam{T_t}_{t\in\R_+}$ of completely positive maps on a \nbd{C^*}algebra. Bhat and Skeide \cite{BhSk00} associate with every CP-semigroup $T$ on a unital \nbd{C^*}al\-ge\-bra $\cB$ a product system $E^\odot=\bfam{E_t}_{t\in\R_+}$ of correspondences $E_t$ over $\cB$ and a unit $\xi^\odot=\bfam{\xi_t}_{t\in\R_+}$, such that $T_t=\AB{\xi_t,\bullet\xi_t}$ and such that $E^\odot$ is generated by $\xi^\odot$.

We refer to the pair $(E^\odot,\xi^\odot)$ as the \hl{GNS-construction} for $T$, and to $E^\odot$ as the \hl{GNS-system}. The unit $\xi^\odot$ is unital if and only if $T$ is a \hl{Markov semigroup} (that is, $T_t(\U)=\U$ for all $t\in\R_+$).

If $T$ is uniformly continuous, then, by Theorem \ref{unitCPSthm}, $E^\odot$ comes along with a unique continuous structure making $\xi^\odot$ a continuous section. Actually, the discussion following Theorem \ref{unitCPSthm} implies that strong continuity of $T$ is sufficient.
\eex

\bemp[Fact.~]
Clearly, if $E^\odot$ is spatial, or even if it just embeds into a continuous spatial product system, then, by Observation \ref{spatUCPSob}, every unit $\xi'^\odot\in CS(E^\odot)$ is continuous. In particular, $T$ must be uniformly continuous.

The other way round, the GNS-system of a strongly continuous but not uniformly continuous CP-semigroup, or every continuous product system containing it, is type III. In particular, it cannot be spatial.
\eemp

\bex\label{E0BPSex}
Let $\vt$ be an \hl{\nbd{E_0}semigroup} (that is, a semigroup of unital endomorphisms) on the unital \nbd{C^*}algebra $\cB$. By $\cB_t$ we denote the right Hilbert \nbd{\cB}module $\cB$ with left action via $\vt_t$, that is, $b.x_t:=\vt_t(b)x_t$. The $\cB_t$ form a product system $\cB^\odot=\bfam{\cB_t}_{t\in\R_+}$ via $u_{s,t}\colon x_s\odot y_t\mapsto \vt_t(x_s)y_t$. It is easy to see that every product system of correspondences over $\cB$ that are one-dimensional as right modules, arises in that way from an \nbd{E_0}semigroup on $\cB$. (Recall that correspondences, by definition, have nondegenerate left action.)

Note that $\cB^\odot$ has at least one unit $\U^\odot=\bfam{\U}_{t\in\R_+}$. In fact, $(\cB^\odot,\U^\odot)$ is the GNS-construction for the Markov semigroup $\vt$. The unit $\U^\odot$ is continuous if and only if $\vt$ is uniformly continuous. If $\vt$ is strongly continuous, then $\U^\odot$ induces a unique continuous structure on $\cB^\odot$. If $\vt$ is not uniformly continuous, then $\cB^\odot$ is type III.
\eex

\bemp[Fact.~]
Clearly, every (continuous) unit in an Arveson system is central (and normalizable). So the spatial Arveson systems are precisely the type I and the type II systems.
\eemp

\bemp[Fact.~]
The following example shows that a type I product system need not be spatial. Once we have a nonspatial type I product system and tensor it with a type II Arveson system (in the obvious way as \it{external} tensor product), that tensor product is a nonspatial type II product system. We do not know an explicit example which cannot be obtained in that way.
\eemp

\bemp[{\protect\cite[Example 4.2.4]{BBLS04}}.~]\label{Inonspatex}
Let $\cB=\sK(H)+\C\U\subset\sB(H)$ be the unitization of the compact operators on some infinite-dimensional Hilbert space $H$. Let $h\in\sB(H)$ be a self-adjoint operator and define the uniformly continuous unital automorphism group $\vt_t=e^{ith}\bullet e^{-ith}$ on $\cB$. Let $\cB^\odot$ be the product system according to Example \ref{E0BPSex} and suppose $\om^\odot$ is a central unital unit. A central element $\om_t\in\cB_t$ fulfills
\beqn{
b.\om_t
~=~
e^{ith}be^{-ith}\om_t
~=~
\om_tb
\text{~~~~~~or~~~~~~}
be^{-ith}\om_t
~=~
e^{-ith}\om_tb
}\eeqn
for all $b\in\cB$. In other words, since the center of $\cB$ is trivial, $e^{-ith}\om_t$ is a multiple of the identity so that $\om_t\in\cB_t=\cB$ is a multiple of $e^{ith}$. Since all $\om_t$ are different from $0$, it follows that $e^{ith}\in\cB$ for all $t\in\R_+$. It follows that $h=-i\frac{d(e^{ith})}{dt}\big|_{t=0}$ is in $\cB$, too.

The other way round, for every $h\notin\cB$ the product system $\cB^\odot$ is type I but nonspatial.
\eemp

\bemp[Fact.~]
By Skeide \cite[Theorem 6.3]{Ske06d}, spatial type I product systems consist of time ordered Fock modules (see Section \ref{counter}). This includes, as a special case, Arveson's result \cite{Arv89} that type I Arveson systems consist of symmetric Fock spaces.

On the contrary, subsystems of type I Arveson systems are type I and, therefore, Fock.
\eemp

New in these notes, we contribute the following results to the classification of product systems.

\bemp[Fact.~]
In Section \ref{counter} we provide a subsystem of a spatial type I (and, therefore, Fock) product system, that is not spatial. The subsystem being the product system of a \it{spatial} CP-semigroup, shows that spatial CP-semigroups need not have spatial GNS-systems. In fact, we prove that spatial CP-semigroups are those with a GNS-systems that embeds into a spatial one; Theorem \ref{CPspatPSthm}.
\eemp

\bemp[Fact.~]
Applying the results from Section \ref{spatialCP}, we show that the nonspatial product systems discussed in Example \ref{Inonspatex} are not even subsystems of spatial systems; Theorem \ref{tInspatthm}.
\eemp

\section{The counter example}\label{counter}

Let $F$ be a correspondence over a unital \nbd{C^*}algebra $\cB$. The \hl{full Fock module}
over $L^2(\R_+,F)$, the completion of the space of (right continuous) \nbd{F}valued step functions, is defined as
\beqn{
\cF(L^2(\R_+,F))
~:=~
\bigoplus_{n=0}^\infty L^2(\R_+,F)^{\odot n}
}\eeqn
where $L^2(\R_+,F)^{\odot 0}=\cB$. By $\om:=\U\in L^2(\R_+,F)^{\odot 0}$ we denote the \hl{vacuum}. The \hl{\nbd{n}particle sector} $L^2(\R_+,F)^{\odot n}$ may be considered as the completion of the space of step functions on $\R_+^n$ with values in $F^{\odot n}$.

Let $\Delta_n$ denote the indicator function of the subset
$\{(t_n,\ldots,t_1)\colon t_n>\ldots>t_1\ge0\}$ of $\R_+^n$.
Then $\Delta_n$ acts on the \nbd{n}particle sector as a projection via
pointwise multiplication (and $\Delta_0$ acts as identity on the vacuum). Set
$\Delta:=\bigoplus\limits_{n=0}^\infty\Delta_n$. The \hl{time ordered
Fock module} is the subcorrespondence
\beqn{
\DG(F)
~:=~
\Delta\cF(L^2(\R_+,F))
}\eeqn
of $\cF(L^2(\R_+,F))$. By $\DG_t(F)$ we denote the subcorrespondence of those functions that are zero if the maximum time argument is $t_n\ge t\in\R_+$. Setting
\beqn{
\SB{u_{s,t}(F_s^m\odot G_t^n)}(s_m,\ldots,s_1,t_n,\ldots,t_1)
~:=~
F_s^m(s_m-t,\ldots,s_1-t)\odot G_t^n(t_n,\ldots,t_1),
}\eeqn
we define bilinear unitaries $u_{s,t}\colon\DG_s(F)\odot\DG_t(F)\rightarrow\DG_{s+t}(F)$ that turn the family $\DG^\odot(F)=\bfam{\DG_t(F)}_{t\in\R_+}$ into a product system, the \hl{time ordered product system} over $F$.

It is easy to check that for every element $\zeta$ in $F$, the elements $\xi_t$ that in each \nbd{n}particle sector assume the constant value $\zeta^{\odot n}$, form a unit $\xi^\odot=\bfam{\xi_t}_{t\in\R_+}$. This unit is also continuous. (See Liebscher and Skeide \cite{LiSk01} for the precise form of all continuous units.) For $\zeta=0$ we obtain the vacuum unit $\om^\odot=\bfam{\om_t}_{t\in\R_+}$ with $\om_t=\U$. The vacuum unit is central and unital. Therefore, if we find a subsystem that does not contain any central unit vector, then this subsystem cannot be Fock.

\lf
Let $\cB=C_0[0,\infty)+\C\U$ denote the unital \nbd{C^*}algebra of all continuous functions on $\R_+$ that have a limit at infinity. Define the Hilbert \nbd{\cB}module $F:=\cB$. We turn $F$ into a correspondence over $\cB$ by defining the left action
\beqn{
b.x
~:=~
\s_1(b)x,
}\eeqn
where $\s_t$ is the left shift by $t$, which acts as $\SB{\s_t(b)}(s)=b(s+t)$. Denote by $\xi^\odot$ the unit corresponding to the parameter $\zeta:=\U\in F$.

\bthm\label{NCexthm}
The product subsystem of $\DG^\odot(F)$ generated by the continuous unit $\xi^\odot$ has no central unit vectors. In particular, it is not isomorphic to a time ordered system.
\ethm

\proof
$F^{\odot n}$ is $\cB$ as Hilbert right module but with left action $b.x=\s_n(b)x$. No nonzero element of $F^{\odot n}$ $(n\ge1)$ can commute with all elements of $\cB$. Therefore, for each $t\ge0$ the set of central elements of $\DG_t(F)$ is the vacuum or \nbd{0}particle sector $\cB$. Commutative \nbd{C^*}algebras do not possess proper isometries. So, the only unit vectors in $\cB$ are unitaries. By multiplying (from the right) with the adjoint, we may assume that such a unit vector is $\U$.

The product subsystem of $\DG^\odot(F)$ generated by $\xi^\odot$ is $E^\odot=\bfam{E_t}_{t\ge0}$ with
\beqn{
E_t
~=~
\cls\BCB{b_n\xi_{t_n}\ldots b_1\xi_{t_1}b_0\colon n\in\N,t_i>0,t_1+\ldots+t_n=t,b_i\in\cB}.
}\eeqn
for $t>0$. Denote by $P_0$ and $P_1$ the projection onto the vacuum component and onto the one-particle component, respectively. We are done if we show that if an element in $x_t\in E_t$ has vacuum component $P_0x_t=\U$, then the one-particle component $P_1x_t\in L^2(\RO{0,t},F)$ is nonzero, too.

Any $x_t\in E_t$ can be approximated by expressions of the form
\beqn{
X_t
~=~
\sum_{i=1}^mb_{n(i)}^{(i)}\xi_{t_{n(i)}^{(i)}}\ldots b_1^{(i)}\xi_{t_1^{(i)}}b_0^{(i)}.
}\eeqn
For $\ve>0$ suppose that $\snorm{x_t-X_t}\le\ve$. Therefore, also $\snorm{P_0x_t-P_0X_t}\le\ve$ and $\snorm{P_1x_t-P_1X_t}\le\ve$. Further, suppose  that $P_0x_t=\U$, that is, suppose that
\beqn{
\Bnorm{\U-\sum_{i=1}^mb_{n(i)}^{(i)}\ldots b_1^{(i)}b_0^{(i)}}
~\le~
\ve.
}\eeqn
The one-particle component of an expression like $b_n\xi_{t_n}\ldots b_1\xi_{t_1}b_0$ is the same as the one-particle component of
\beqn{
b_n(\U\oplus\I_{\RO{0,t_n}}\U)\ldots b_1(\U\oplus\I_{\RO{0,t_1}}\U)b_0,
}\eeqn
where $\I_A$ denotes the \hl{indicator function} of the set $A$. The one-particle component of this expression is
\beqn{
\I_{\RO{t_1+\ldots t_{n-1},t_1+\ldots t_n}}\s_1(b_n)b_{n-1}\ldots b_1b_0
~+~
\ldots
~+~
\I_{\RO{t_1+t_2,t_1}}\s_1(b_n\ldots b_2)b_1b_0
~+~
\I_{\RO{t_1,0}}\s_1(b_n\ldots b_1)b_0.
}\eeqn
From $\lim_{s\to\infty}\SB{\s_1(b)c}(s)=\lim_{s\to\infty}b(s+1)c(s)=\lim_{s\to\infty}b(s)c(s)=\lim_{s\to\infty}\SB{bc}(s)$ and the fact that $X_t$ contains only finitely many summands it follows that
\beqn{
\lim_{s\to\infty}\AB{P_1X_t,P_1X_t}(s)
~=~
\lim_{s\to\infty}t\Babs{\sum_{i=1}^mb_{n(i)}^{(i)}\ldots b_1^{(i)}b_0^{(i)}}^2\!\!(s)
~=~
\lim_{s\to\infty}t\AB{P_0X_t,P_0X_t}(s).
}\eeqn
The function $\AB{P_0X_t,P_0X_t}$ of $s$ is uniformly close to $1$. So, $\snorm{P_1X_t}^2\ge\lim_{s\to\infty}\AB{P_1X_t,P_1X_t}(s)$. Therefore, $\norm{P_1x_t}$ is bigger than a number arbitrarily close to $t\ne0$.\qed

\section{Spatial CP-semigroups}\label{spatialCP}

Recall that a CP-map $T$ \hl{dominates} another $S$, if the difference $T-S$ is a CP-map, too. A CP-semigroup $T$ \hl{dominates} another $S$, if $T_t$ dominates $S_t$ for all $t\in\R_+$. A CP-semigroup $S$ on a unital \nbd{C^*}algebra $\cB$ is \hl{elementary}, if it has the form $S_t(b)=c_t^*bc_t$  for some semigroup $c=\bfam{c_t}_{t\in\R_+}$ of elements $c_t$ in $\cB$. Following Arveson's definition in \cite{Arv97a} for $\cB=\sB(H)$, we say a \hl{unit} for $T$ is a semigroup $c$ in $\cB$ such that $T$ dominates the elementary CP-semigroup $b\mapsto c_t^*bc_t$.

Arveson defines a CP-semigroup to be \it{spatial}, if it admits units. Without continuity conditions, every CP-semigroup dominates an elementary CP-semigroup, namely, the \nbd{0}semigroup which is $\id_\cB$ for $t=0$ and $0$ otherwise. Depending on the context, there are several topologies  around in which a CP-semigroup can be continuous with respect to time $t\in\R_+$. The uniform (or norm) topology, the strong and weak topologies of operators on the Banach space $\cB$, and pointwise versions of all the operator topologies when $\cB\subset\sB(H)$ is a concrete operator algebra, for instance, if $\cB$ is a von Neumann algebra.

In Arveson's definition, a unit for a CP-semigroup on $\sB(H)$ is required pointwise continuous in the strong operator topology of $\sB(H)$. The usual topology used for CP-semi\-groups on a \nbd{C^*}algebra is the strong topology of operators on the Banach space $\cB$. It is well-known that a weakly continuous semigroup is also strongly continuous; see, for instance, Bratteli and Robinson \cite[Corollary 3.1.8]{BrRo87}. For semigroups $c=\bfam{c_t}_{t\in\R_+}$ in $\cB$, in absence of a strong topology on $\cB$ or a predual $\cB_*$, the only obvious topology apart from the norm topology is the weak topology. However, if $t\mapsto c_t$ is weakly continuous, then also the semigroup $\bullet c_t\colon b\mapsto bc_t$ of operators on $\cB$ is weakly, hence, strongly continuous. However, strong continuity for that semigroup means that, in particular for $b=\U$, the map $t\mapsto c_t$ is norm continuous.

\brem
This is not a contradiction to the existence of weakly continuous unitary groups on a Hilbert space $H$. Here weakly continuous refers to the weak operator topology of $\sB(H)$, which is much weaker than the weak topology of $\sB(H)$.
\erem

\brem
It should be noted that, apart from Arveson's, there is another definition of \it{spatial} CP-semigroups on $\sB(H)$ due to Powers \cite{Pow04}. Powers' definition is much more restrictive, and only Arveson's definition gives equivalence of the notions of spatiality for product systems of Hilbert spaces and for CP-semigroups on $\sB(H)$.
\erem

\bdefi\label{spatCPdefi}
A strongly continuous CP-semigroup on a unital \nbd{C^*}algebra is \hl{spatial}, if it admits a continuous unit.
\edefi

To prove the following theorem, we have to recall a few definitions and facts from Barreto, Bhat, Liebscher, and Skeide \cite{BBLS04}. A \hl{kernel} over a set $S$ with values in the set $\sB(\cA,\cB)$ of bounded mappings from $\cA$ to $\cB$ is just a map $\eK\colon S\times S\rightarrow\sB(\cA,\cB)$. The kernel $\eK$ is \hl{completely positive definite} (or \hl{CPD}), if
\beqn{
\sum_{i,j}b_i^*\eK^{s_i,s_j}(a_i^*a_j)b_j
~\ge~
0
}\eeqn
for all choices of finitely many $s_i\in S,a_i\in\cA,b_i\in\cB$. If $\cA=\cB$, we say $\eK$ is a kernel \hl{on} $\cB$. A \hl{CPD-semigroup} on $\cB$ is a family $\eT=\bfam{\eT_t}_{t\in\R_+}$ of CPD-kernels on $\cB$, such that for all $s,s'\in S$ the maps $\eT_t^{s,s'}$ form semigroups on $\cB$. The CPD-semigroup is \hl{continuous} in a certain topology, if every semigroup $\eT_t^{s,s'}$ is continuous in that topology.

If $E^\odot$ is a product system and $S$ a set of units, then it is easy to check that the semigroups $\eT^{\xi,\xi'}=\AB{\xi_t,\bullet\xi'_t}$ (see Section \ref{prel}) form a CPD-semigroup over $S$. Generalizing the result for CP-semigroups discussed in Example \ref{CPGNSex}, the result \cite[Theorem 4.3.5]{BBLS04} asserts that for every CPD-semigroup on a unital \nbd{C^*}algebra $\cB$ there is a product system $E^\odot$ and a family $\bfam{{\xi^s}^\odot}_{s\in S}$ of units for $E^\odot$, such that
\beqn{
\eT_t^{s,s'}
~=~
\AB{\xi^s_t,\bullet\xi^{s'}_t}
}\eeqn
for all $s,s'\in S$. The subsystem of $E^\odot$ generated by all these units is unique in an obvious way, and we refer to it as the \hl{GNS-system} of $\eT$. Once again, as in Example \ref{CPGNSex}, if $\eT$ is strongly continuous, then the GNS-system of $\eT$ is continuous.

\bthm\label{CPspatPSthm}
A strongly continuous CP-semigroup is spatial, if and only if its GNS-system can be embedded into a continuous spatial product system.
\ethm

\proof
``$\Longleftarrow$.''~
Let $T$ be a strongly continuous CP-semigroup. Suppose $E^\odot$ is a continuous product system with a unit $\xi^\odot\in CS(E^\odot)$ such that $T_t=\AB{\xi_t,\bullet\xi_t}$ and a central unital unit $\om^\odot\in CS(E^\odot)$. Then $c_t:=\AB{\om_t,\xi_t}$ is a semigroup of elements in $\cB$; see \cite[Section 5.1]{BBLS04}. As explained in the beginning of this section, since $t\mapsto\AB{\om_t,\bullet\xi_t}=\bullet c_t$ is strongly continuous, $t\mapsto c_t$ is uniformly continuous. Define the bilinear projection $q_t:=\id_t-\om_t\om_t^*\in\cB^{a,bil}(E_t)$. By
\beqn{
T_t(b)-c_t^*bc_t
~=~
\AB{\xi_t,b\xi_t}-\AB{\xi_t,\om_t}b\AB{\om_t,\xi_t}
~=~
\AB{\xi_t,q_tb\xi_t}
~=~
\AB{(q_t\xi_t),b(q_t\xi_t)},
}\eeqn
we see that $T_t-c_t^*\bullet c_t$ is completely positive for all $t\in\R_+$.

``$\Longrightarrow$.''~
Let $c=\bfam{c_t}_{t\in\R_+}$ be a unit for the strongly continuous CP-semigroup $T$. Then the strongly continuous semigroup $\eT_t$ of kernels over $\CB{0,1}$ on $\cB$ defined by setting
\beqn{
\Matrix{\eT^{0,0}_t&\eT^{0,1}_t\\\eT^{1,0}_t&\eT^{1,1}_t}
~:=~
\Matrix{\id_\cB&\bullet c_t\\c_t^*\bullet&T_t}
~:=~
\Matrix{0&0\\0&T_t-c_t^*\bullet c_t}+\Matrix{\id_\cB&\bullet c_t\\c_t^*\bullet&c_t^*\bullet c_t}
}\eeqn
is completely positive definite. (Indeed, $c$ is a unit for $T$. So, the first summand is CPD. The second is a simple example of what is called an \hl{elementary} CPD-semigroup in Skeide \cite{Ske08p4}. It is, clearly, CPD.) The GNS-system of $\eT$ is, then, a spatial continuous product system with unital  central unit ${\xi^0}^\odot$ containing the GNS-system of $T_t$ as the subsystem generated by the unit ${\xi^1}^\odot$.\qed

\brem
Theorem \ref{NCexthm} tells us that we may \bf{not} replace Definition \ref{spatCPdefi} with the property that $T$ has a spatial GNS-system. Theorem \ref{CPspatPSthm} tells us that we \bf{may} replace Definition \ref{spatCPdefi} with the property that the GNS-system of $T$ embeds into a spatial product system. The clarification of these facts was among the main scopes of these notes.
\erem

Directly from Observation \ref{spatUCPSob} we conclude:

\bcor\label{spatembcor}
Every spatial strongly continuous CP-semigroup is uniformly continuous. The other way round, every strongly continuous CP-semigroup that is not uniformly continuous, is nonspatial, too.
\ecor

The opposite statement need not be true. In fact, the following sharp version of the corollary allows us to show in Theorem \ref{tInspatthm} that Example \ref{E0BPSex} is a counter example.

Recall that the generator $L:=\frac{d}{dt}\big|_{t=0}T_t$ of a uniformly continuous CP-semigroup on $\cB$ has \hl{Christensen-Evans form}, if there are a correspondence $F$ over $\cB$, an element $\zeta\in F$, and an element $\beta\in\cB$, such that $L(b)=\AB{\zeta,b\zeta}+b\beta+\beta^*b$. We also say $L$ is a \hl{Christensen-Evans generator}.

\bcor\label{spatCEcor}
A strongly continuous CP-semigroup on a unital \nbd{C^*}algebra is spatial if and only if it has a Christensen-Evans generator.
\ecor

\proof
If the semigroup is spatial, then by Corollary \ref{spatembcor} it is uniformly continuous. Therefore, by \cite[Lemma 5.1.1]{BBLS04} its generator has Christensen-Evans form. On the other hand, by \cite[Corollary 5.1.3]{BBLS04}, the product system of a CP-semigroup with Christensen-Evans generator embeds into a time ordered product system, so that, by Theorem \ref{CPspatPSthm}, the CP-semigroup is spatial.\qed

\lf
As a simple corollary, we derive existence of a type I product system that does not embed into a time ordered system.

\bthm\label{tInspatthm}
Like in Example \ref{Inonspatex}, let $\cB=\sK(H)+\C\U\subset\sB(H)$ and consider the automorphism semigroup $\vt_t=e^{ith}\bullet e^{-ith}$ for a self-adjoint operator $h\in\sB(H)\backslash\cB$.

Then $\vt$ is a uniformly continuous Markov semigroup that is not spatial. Equivalently, its GNS-system, though type I, does not embed into a time ordered system.
\ethm

\proof
$\vt$ has the generator $L(b)=i(hb-bh)$. Let us suppose that $L$ can be written in the Christensen-Evans form $L(b)=\AB{\zeta,b\zeta}+b\beta+\beta^*b$. Calculating $\sum_{i,j}b_i^*L(a_i^*a_j)b_j$ in either way, one easily verifies that $\sum_ia_ib_i=0$ $\Longrightarrow$ $\sum_ia_i\zeta b_i=0$ for all finite choices of $a_i,b_i\in\cB$. In particular, $b\zeta=\zeta b$, so that $\zeta$ and $\AB{\zeta,\zeta}$ are central. Since $\vt$ is Markov, $\Re\beta=-\frac{1}{2}\AB{\zeta,\zeta}$. Consequently, the Christensen-Evans form of $L$ simplifies to $L(b)=i(h'b-bh')$, where $h':=-\Im\beta\in\cB$.

Exponentiating the two forms of $L$, we find $e^{ith}be^{-ith}=e^{ith'}be^{-ith'}$ or $e^{-ith'}e^{ith}b=be^{-ith'}e^{ith}$ for all $b\in\cB$. In other words, the unitary $e^{-ith'}e^{ith}\in\sB(H)$ must be a multiple of the identity. It follows $h=h'+\vp\U$ for some real number $\vp$. In particular, $h\in\cB$.

In other words, if $h\notin\cB$, then $L$ cannot be written in Christensen-Evans form, so that $\vt$ is nonspatial.\qed

\brem
Following Definition \ref{spatCPdefi}, Corollaries \ref{spatembcor} and \ref{spatCEcor} are intrinsic statements about strongly continuous CP-semigroups. However, we think it would not be easily possible to prove these statements without reference to product systems. These results continue a whole series of intrinsic statements about CP- or CPD-semigroups that have comparably simple proofs in terms of their GNS-system; see also Liebscher and Skeide \cite[Remark 3.6]{LiSk08}.
\erem

\newpage

\section{Remarks}\label{rem}

\brem
After the short introduction, in Section \ref{prel} we went straight to our definition of product systems, in particular spatial ones. Of course, this ignores historical order. But, giving an account on the historic development would require several pages and is not among the scopes of these notes.

Furthermore, in Example \ref{CPGNSex} we discussed the relation of a CP-semigroup and its GNS-system. Also here we did not discuss the Arveson system of a CP-semigroup on $\sB(H)$ defined by Bhat \cite{Bha96}. This construction of the Arveson system of a CP-semigroup on $\sB(H)$ is rather indirect. (It is the product system of the so-called \it{minimal dilation} of the CP-semigroup.) Additionally, it is \bf{not} the product system we mean, when we speak about the GNS-system of a CP-semigroup. (The latter is a product system of correspondences over $\sB(H)$, while the former is an Arveson system.) Depending on the construction of the Arveson system of the minimal dilation, the relation between the GNS-system and the Arveson system of a CP-semigroup on $\sB(H)$ is either by \it{Morita equivalence} of product systems or by the \it{commutant} of product systems; see \cite{Ske09,Ske08p1}. We do not have enough space to explain all these peculiarities, and they are beyond what is needed for these notes.
\erem

\brem
Product systems, in their concise formulation as a sort of ``measurable semigroup under tensor product'' of Hilbert spaces, occurred first in Arveson \cite{Arv89} in the study of \nbd{E_0}semigroups on $\sB(H)$. Like Arveson's definition, the definition of continuous product systems in \cite{Ske03b} is motivated by the relation between \nbd{E_0}semigroups and product systems. The point in \cite{Ske03b} is that a product system of a strongly continuous \nbd{E_0}semigroup is continuous. The validity of this definition has been confirmed in \cite{Ske07,Ske08p1,Ske09a}, showing that every (full) product system comes from an \nbd{E_0}semigroup, and that (under countability hypothesis) \nbd{E_0}semigroups are classified up to \it{stable} cocycle conjugacy. This is in full correspondence with Arveson's theory \cite{Arv89,Arv90a,Arv89a,Arv90} for Arveson systems.

If a product system has a unital unit, then it is easy to construct an \nbd{E_0}semigroup for that product system. It should be noted that all proofs of the statements in Section \ref{prel} about continuous product systems in the presence of a continuous unital unit make use of that construction.
\erem

\brem
Note that product systems can easily be defined also for nonunital $\cB$. Recently, Skeide \cite{Ske09a} generalized the relation between (continuous) product systems and (strongly continuous) \nbd{E_0}semigroups to the case of \nbd{\sigma}unital $\cB$. However, units have been defined, so far, only for unital $\cB$ with the stated marginal condition $\xi_0=\U$. We are convinced that in a reasonable definition of unit for nonunital $\cB$, the unit must consist of maps in $\sB^a(\cB,E_t)$, the \it{strict} completion of $E_t$; see, for instance, the discussion in \cite[Section 7]{Ske09}. A systematic extension of the theory in that direction is largely missing.
\erem

\brem
Many counter examples in the \nbd{C^*} case lose their validity for the von Neumann case, and the classification becomes again simpler. The types of product systems are defined as for \nbd{C^*}cor\-re\-spond\-ences, that is, referring to continuous sets of units. The only difference is that now \it{generating} means generating in the strong topology. We give brief account including old results, mainly from \cite{BBLS04}, but also results in \cite{Ske08p1} obtained a quite a while after the first version of these notes.

The main result of \cite{BBLS04} can be phrased as follows: Non-type III systems of von Neumann correspondences are spatial. In particular, type I systems are Fock. The proof in \cite{BBLS04} goes by establishing equivalence with the main result of Christensen and Evans \cite{ChrEv79}: Bounded derivations with values in a von Neumann correspondence are inner.

It is an open question, if a subsystem of a type I system must be type I, too. Likewise, for type II systems.

With the first concise definition of \it{strongly continuous} product systems of von Neumann correspondences, in \cite{Ske08p1} a good deal of the relevant statements of the \nbd{C^*}theory could be put through for von Neumann correspondences. (Note that von Neumann modules are modules of operators, and strongly continuous, here, refers to the strong operator topology.)

Theorem \ref{unitCPSthm} remains true for strongly continuous product systems and strongly continuous units. This means, the just stated results from \cite{BBLS04} remain true also for strongly continuous product systems: Non-type III systems are spatial and type I systems are Fock.

In \cite{Ske08p1}, where also the case of spatial CP-semi\-groups on a von Neumann algebra is treated, the following is proved: A strongly continuous normal CP-semigroup on a von Neumann algebra is spatial if and only if its GNS-system (of von Neumann correspondences) is spatial. This is in contrast with the counter example for the \nbd{C^*}case in Section \ref{counter}.

However, we know from Fagnola, Liebscher, and Skeide \cite{FLS05p,Ske04} that the GNS-systems of the Markov semigroup of Brownian motion or of the Ornstein-Uhlenbeck processes are nonspatial. Therefore, \it{strong type I} product systems (in the sense that they are generated by a strongly operator continuous set of units) need not be spatial.

Under strong closure, the semigroup in Example \ref{Inonspatex} becomes an inner automorphism semigroup. Therefore, the product system is the trivial one and, in particular, Fock (over the one-particle sector $\zero$).  This is in correspondence with the result that, for von Neumann correspondences, type I implies Fock.
\erem

\brem
In Corollary \ref{spatembcor} we have seen that spatial strongly continuous CP-semigroups on a unital \nbd{C^*}algebra $\cB$ are uniformly continuous. This is so due to the fact that there is no semigroup $c=\bfam{c_t}_{t\in\R_+}$ in $\cB$ continuous in any of the natural topologies of $\cB$, that would not be uniformly continuous. Also for a von Neumann algebra $\cB$, from the beginning, it is not reasonable to consider CP-semigroups that are strongly continuous. A result due to Elliott \cite{Ell00} asserts that such a semigroup would be uniformly continuous. But in weaker topologies where also units $c^*$ need no longer be uniformly continuous, there are much richer classes of spatial CP-semigroups. Actually, practically all known explicit examples of CP-semigroup on $\sB(H)$ are spatial in this sense, when continuity is with respect to the strong operator topology.
\erem

\brem
Several more results have been obtained. Spatial Markov semigroups ($C^*$ or von Neumann case) are precisely those that admit \it{Hudson-Parthasarathy} dilations, that is, dilations by cocycle perturbations of \it{noises}; Skeide \cite[Sections 10 and 13]{Ske08p1}. The \it{Powers sum} \cite{Pow04} of CP-semigroups that are spatial in the (restrictive) sense of Powers, has been generalized not only to all spatial CP-semigroups (in our sense) on $\sB^a(E)$ but also to all \it{spatial} CPD-semigroups on $\cB$; see Skeide \cite{Ske08p4}. It is also shown that the \it{spatialized} GNS-system of the generalized Powers sum is always the \it{product} (in the sense of \cite{Ske06d}) of the spatialized GNS-systems of the summands.
\erem

\setlength{\baselineskip}{2.5ex}

\newcommand{\Swap}[2]{#2#1}\newcommand{\Sort}[1]{}
\providecommand{\bysame}{\leavevmode\hbox to3em{\hrulefill}\thinspace}
\providecommand{\MR}{\relax\ifhmode\unskip\space\fi MR }
\providecommand{\MRhref}[2]{%
  \href{http://www.ams.org/mathscinet-getitem?mr=#1}{#2}
}
\providecommand{\href}[2]{#2}

\noindent
B.V.\ Rajarama Bhat: \it{Statistics and Mathematics Unit, Indian Statistical Institute Bangalore, R.\ V.\ College Post, Bangalore 560059, India},
E-mail: \tt{bhat@isibang.ac.in},\\
Homepage: \tt{http://www.isibang.ac.in/Smubang/BHAT/}

\lf\noindent
Volkmar Liebscher: \it{Institut für Mathematik und Informatik, Ernst-Moritz-Arndt-Universität Greifswald, 17487 Greifswald, Germany},\\
E-mail: \tt{volkmar.liebscher@uni-greifswald.de},\\
Homepage: \tt{http://www.math-inf.uni-greifswald.de/biomathematik/liebscher/}

\lf\noindent
Michael Skeide: \it{Dipartimento S.E.G.e S., Università degli Studi del Molise, Via de Sanctis, 86100 Campobasso, Italy},
E-mail: \tt{skeide@math.tu-cottbus.de},\\
Homepage: \tt{http://www.math.tu-cottbus.de/INSTITUT/lswas/\_skeide.html}

\end{document}